\documentclass[a4paper,10pt,fleqn,conference]{ieeeconf}
\IEEEoverridecommandlockouts
\usepackage{graphicx}
\usepackage{macros_conf}
\usepackage{amsfonts}
\usepackage{nccmath}
\usepackage{amsmath}
\usepackage{amssymb}
\usepackage{mathtools}
\usepackage{arydshln}
\usepackage{here}
\usepackage{macros_e}
\usepackage{mathfonts}
\usepackage{footnote}

\input{colordvi}

\everymath={\displaystyle}
\newcommand{\argmax}{\mathop{\rm argmax}\limits}

\title{
\LARGE \bf 
Local Lipschitz Constant Computation of ReLU-FNNs:\\
Upper Bound Computation with Exactness Verification
}
\author{Yoshio Ebihara, Xin Dai, Victor Magron, Dimitri Peaucelle, Sophie Tarbouriech
\thanks{
Y. Ebihara and X. Dai are with the 
Graduate School of Information Science and 
Electrical Engineering, Kyushu University, 
744 Motooka, Nishi-ku, Fukuoka 819-0395, Japan, 
Y. Ebihara was also with 
LAAS-CNRS, Universit\'{e} de Toulouse, CNRS, Toulouse, France,   
in 2011.  
V. Magron, D. Peaucelle, and S. Tarbouriech are
LAAS-CNRS, Universit\'{e} de Toulouse, CNRS, Toulouse, France.  
This work was supported by JSPS KAKENHI Grant Number JP21H01354.  
This work was also supported by the AI Interdisciplinary Institute ANITI funding, through the French "Investingfor the Future PIA3" program under the Grant agreement n°ANR-19-PI3A-0004 as well as by the National Research Foundation, Prime Minister's Oﬃce, Singapore under its Campus for Research Excellence and Technological Enterprise (CREATE) programme. 
}%
}

\begin{document}
\maketitle
\thispagestyle{empty}
\pagestyle{empty}

\begin{abstract}
This paper is concerned with the computation of 
the local Lipschitz constant of 
feedforward neural networks (FNNs) with activation functions being rectified linear units (ReLUs).  
The local Lipschitz constant of an FNN for a target input 
is a reasonable measure for its quantitative evaluation of the reliability.  
By following a standard procedure using multipliers that capture the behavior of ReLUs,
we first reduce the upper bound computation problem of the local Lipschitz constant
into a semidefinite programming problem (SDP).  
Here we newly introduce copositive multipliers to capture the ReLU behavior accurately.  
Then, by considering the dual of the SDP for the upper bound computation, 
we second derive a viable test to conclude the exactness of the computed upper bound.  
However, these SDPs are intractable for practical FNNs with hundreds of ReLUs.  
To address this issue, we further propose a method to construct a reduced order model 
whose input-output property is identical to the original FNN over a neighborhood of the target input.  
We finally illustrate the effectiveness of the model reduction and exactness verification methods 
with numerical examples of practical FNNs.  

\end{abstract}

\section{Introduction}
\label{sec:intro}

Control theoretic approaches for the reliability certification of 
deep neural networks have gained considerable attention recently.  
These studies are roughly classified into two categories: 
the treatment of static feedforward neural networks (FNNs) 
\cite{Raghunathan_NIPS2018,Raghunathan_ICLR2018,Fazlyab_arxiv2019,Fazlyab_IEEE2022,Chen_NIPS2020} 
and the treatment of dynamical networks such as recurrent neural networks (RNNs)
\cite{Yin_IEEE2022,Scherer_IEEEMag2022,Revay_LCSS2021,Ebihara_EJC2021}.  
In the present paper, we focus on FNNs, and investigate 
the computation problem of the local Lipschitz constant for a target input.  
The effectiveness of FNNs in image recognition and pattern classification, etc.,  is widely recognized.  
However, it is known that some FNNs
exhibit unreliable behavior such that small perturbation
on an input leads to a totally different output \cite{Eykholt_arxiv}. 
The existence of such inputs, known as adversarial inputs (perturbations), 
shows the unreliability of the FNN.  
By computing the local Lipschitz constant, 
we can ensure the absence of such adversarial inputs, 
and thereby certify the reliability of the FNN.  

In this paper, we first reduce the upper bound computation problem
of the local Lipschitz constant into a semidefinite programming problem (SDP), 
by following the standard procedure using multipliers that capture
the input-output behavior of ReLUs in quadratic constraints \cite{Fazlyab_IEEE2022}.  
Here we newly introduce a set of multipliers constructed from copositive matrices \cite{Dur_2010}, 
in order to accurately capture the nonnegative behavior of ReLUs.  
In addition, we prove that the copositive multipliers encompass existing multipliers
such as O'Shea-Zames-Falb multipliers 
\cite{O'Shea_IEEE1967,Zames_SIAM1968,Carrasco_EJC2016,Fetzler_IFAC2017}
and the multipliers introduced in \cite{Fazlyab_IEEE2022}.    
We next focus on the dual of the SDP for the upper bound computation.  
Then we derive a rank condition on the dual optimal solution under which 
we can conclude the exactness of the obtained upper bound.  
In particular, this also enables us to extract the worst case input 
that deviates the corresponding output most far away from the original one.  
However, for practical FNNs with hundreds of ReLUs, 
these SDPs are inherently intractable, since the size of 
the multiplier grows linearly with the number of ReLUs
and thus the computational complexity for solving these SDPs increases quite rapidly.  
To address this issue, we further propose a method to construct a reduced order model 
whose input-output property is identical to the original FNN over a neighborhood of the target input.  
We finally illustrate the effectiveness of the model reduction and exactness verification methods 
by numerical examples on academic and practical FNNs, 
where in the latter case we use the FNN in \cite{Raghunathan_ICLR2018} 
designed as an MNIST classifier.  

\noindent
{\bf Notation:} 
The set of $n\times m$ real matrices is denoted by $\bbR^{n\times m}$, and
the set of $n\times m$ entrywise nonnegative matrices is denoted
by $\bbR_+^{n\times m}$. 
For a matrix $A$, we also write $A\geq 0$ to denote that  
$A$ is entrywise nonnegative. 
We denote the set of $n\times n$ real symmetric matrices by $\bbS^n$.
For $A\in\bbS^n$, we write $A\succ 0\ (A\prec 0)$ to
denote that $A$ is positive (negative) definite.  
For $A\in\bbR^{n\times n}$, we define $\He\{A\}:=A+A^T$.
For $A\in\bbR^{n\times n}$ and $B\in\bbR^{n\times m}$,
$(\ast)^TAB$ is a shorthand notation of $B^TAB$.  
We denote by $\bbD^n\subset\bbR^{n\times n}$ the set of diagonal matrices.  
A matrix $M\in\bbR^{n\times n}$ is said to be 
Z-matrix if $M_{ij}\le 0$ for all $i\neq j$.  
Moreover, $M$ is said to be doubly hyperdominant if
it is a Z-matrix and 
$M \one_n \ge 0$, $\one_n^T M \ge 0$,   
where $\one_n\in\bbR^n$ stands for the all-ones-vector.  
We denote by $\DHD^n\subset \bbR^{n\times n}$ the set of 
doubly hyperdominant matrices.  
For $v\in\bbR^{n}$, we define $|v|_2:=\sqrt{\sum_{j=1}^{n} v_j^2}$.  
For $v_0\in\bbR^{n}$ and $\varepsilon>0$, we define
$\clB_{\varepsilon}(v_0):=\left\{v\in\bbR^{n}:\ |v-v_0|_2\le \varepsilon \right\}$.  
Finally, for the $i$-th unit vector $e_i\in\bbR^n$ and the index set 
$\clN\subset\left\{1,\cdots,n\right\}$, 
we denote by $\bigoplus_{i \in \clN}e_i\in\bbR^{|\clN|\times n}$
the matrix that arrays $e_i^T$ vertically for $i\in\clN$.  
For instance, if $n=4$ and $\clN = \{1,2\}$ then we have
$
\bigoplus_{i \in \clN}e_i=
\begin{bmatrix}
1 & 0 & 0 & 0 \\
0 & 1 & 0 & 0 \\
\end{bmatrix}
$.     
%

\section{Local Lipschitz Constant Computation Problem of ReLU-FNNs and Basic Results}
\label{sec:RNN}

\subsection{Local Lipschitz Constant Computation Problem}
\label{sub:prob}

Let us consider the single-layer FNN described by
\begin{equation}
G:\ z = \Wout\Phi(\Win w + \bin)+\bout.  
\label{eq:FNN}
\end{equation}
Here $w\in\bbR^m$ and $z\in\bbR^l$ are 
the input and the output of the FNN, respectively.  
On the other hand, $\Wout\in\bbR^{l\times n}$ and
$\Win\in\bbR^{n\times m}$ are constant matrices
constructed from the weightings of the edges in the FNN.  
The constant vectors 
$\bin\in\bbR^n$ and $\bout\in\bbR^l$ are biases at input and output, 
respectively.   
Finally, $\Phi:\ \bbR^n\to \bbR^n$ is the
static activation function typically being nonlinear.  
In the following, we denote by $z=G(w)$ the input-output map of the FNN.  

In this paper, we consider the typical case where the activation function is
the (entrywise)
rectified linear unit (ReLU) whose input-output property is given by
\begin{equation}
\begin{array}{@{}l}
 \Phi(q)=\left[\ \phi(q_1)\ \cdots\ \phi(q_n)\ \right]^T,\\
 \phi: \bbR\to \bbR,\quad
 \phi(\eta)=
  \left\{
   \begin{array}{cc}
    \eta & (\eta\ge 0),  \\
    0 & (\eta< 0).  \\
   \end{array}
  \right.    
\end{array}  
\label{eq:ReLU}
\end{equation}
For $p,q\in\bbR^n$, it is shown in \cite{Raghunathan_NIPS2018,Groff_CDC2019} that 
$p=\Phi(q)$ holds if and only if
\begin{equation}
p-q\ge 0,\ p\ge 0,\ (p-q)\otimes p =0  
\label{eq:ReLU_alg}
\end{equation}
where $\otimes$ stands for the Hadamard product.  

The local Lipschitz constant for the FNN given by \rec{eq:FNN}
at a target input is defined as follows.  

\begin{definition} (Local Lipschitz Constant at Target Input) 
For given target input $w_0\in\bbR^m$ and error bound $\varepsilon>0$, 
the local Lipschitz constant $L_{w_0,\varepsilon}$ 
for the FNN $G$ given by \rec{eq:FNN} and \rec{eq:ReLU}
is defined by
 \begin{equation}
\scalebox{0.92}{$
\begin{array}{@{}l}
L_{w_0,\varepsilon}:= \min\left\{L:\ |G(w)-G(w_0)|_2 \le L\  \forall w\in \clB_{\varepsilon}(w_0) \right\}  
\end{array}$}.  
\label{eq:Lipschitz2}
 \end{equation}
\label{def:Lipschitz}
\end{definition}

The target input $w_0\in\bbR^m$ and the error bound $\varepsilon>0$ 
are determined by the one who evaluates the reliability of the FNN.  
The goal of this paper is to compute 
the local Lipschitz constant as accurately as possible.
From \rdef{def:Lipschitz}, we see that the bias $\bout$ at output
is irrelevant to the local Lipschitz constant and hence we make the following assumption.
\begin{assumption}
We assume $\bout=0$ in \rec{eq:FNN} in the following.  
\end{assumption}

The computation of the local Lipschitz constant is strongly
motivated  by the demand of robustness analysis of
FNNs against adversarial inputs (perturbations).
To be more concrete, let us follow the convincing discussion in \cite{Fazlyab_arxiv2019}
and consider the case where the FNN $G$ given by \rec{eq:FNN} is used
as a classifier.  Namely, the FNN $G$ receives $m$ features as input and
returns $l$ scores as output.
We define the classification rule $C: \bbR^{m}\to \{1,\cdots,l\}$ by
\begin{equation}
C(w):=\argmax_{1\le i \le l} G_i(w).  
\label{eq:classifier}
\end{equation}
Now suppose $w_0$ is an input that is classified correctly
by the classifier.  
Then, the classifier is locally robust at the target input $w_0$
against all the perturbed inputs $w\in\clB_{\varepsilon}(w_0)$ if
\begin{equation}
 C(w)=C(w_0)\ \forall w\in \clB_{\varepsilon}(w_0).   
\label{eq:robl2}
\end{equation}
Regarding this condition for ensuring the robustness of the FNN, 
it has been shown in \cite{Fazlyab_arxiv2019} that 
the next result holds.  
\begin{proposition}\cite{Fazlyab_arxiv2019} 
For the classifier $C: \bbR^{m}\to \{1,\cdots,l\}$ given by
\rec{eq:classifier} and \rec{eq:FNN} and given $w_0\in\bbR^m$ and
 $\varepsilon>0$, let us define $i^\star:=C(w_0)$.
 Then the condition \rec{eq:robl2} holds if
\begin{equation}
 L_{w_0,\varepsilon}\le
 \dfrac{1}{\sqrt{2}}\min_{1\le j\le l, j\ne i^{\star}}
 G_{i^\star}(w_0)-G_j(w_0).  
\label{eq:robl2cond}
\end{equation}
\label{pr:rob}
\end{proposition}

From this proposition, we see that we can examine 
the absence of adversarial inputs and hence ensure the reliability of the FNN
by computing the local Lipschitz constant.

\subsection{Basic Results}
\label{sub:key}

To grasp the properties of the ReLU given by \rec{eq:ReLU}, 
we borrow the idea of 
integral quadratic constraint (IQC) theory \cite{Megretski_IEEE1997}.  
Namely, we introduce the set of multipliers 
$\bPi^\star\subset\bbS^{2n+1}$ such that
\begin{equation}
\scalebox{0.78}{$
\begin{array}{@{}l}
\bPi^\star:=
\left\{
\Pi\in\bbS^{2n+1}:\ 
 \left[
\arraycolsep=0.3mm
 \begin{array}{c}
 1 \\
 q \\
 p \\
 \end{array}
 \right]^T\!\!\!
 \Pi
 \left[
\arraycolsep=0.5mm
 \begin{array}{c}
 1 \\
 q \\
 p \\
 \end{array}
 \right]\ge 0\ \forall q,p\in\bbR^n\ \mbox{s.t.}\ p=\Phi(q)
\right\}.  
\end{array}$}
\label{eq:IQC0}
\end{equation}
This type of multiplier is already introduced in \cite{Fazlyab_arxiv2019,Fazlyab_IEEE2022}.  
Then, we can readily obtain the next results.  
\begin{theorem}
For given input $w_0\in\bbR^m$ and $\varepsilon>0$, 
let us define $z_0:=G(w_0)$ where $G$ is the FNN given by \rec{eq:FNN} and \rec{eq:ReLU}.   
Suppose $\Pi\in\bPi^\star$.  
Then, $L_{w_0,\varepsilon}\le \sqrt{\Lsq}$ holds for the FNN $G$ if there exists $\tau\ge 0$ such that
\begin{equation}
\scalebox{0.95}{$
\begin{array}{@{}l}
(*)^T  
\left[
\begin{array}{ccc} 
-\Lsq+\tau\epsilon^2 & 0 & 0\\
\ast & -\tau I_m &0\\
\ast &\ast & I_l
\end{array}
\right]
\left[
\begin{array}{ccc} 
1&0&0\\ 
-w_0 & I_m &0\\ 
-z_0 & 0 & \Wout 
\end{array}\right] \vspace*{1mm}\\
\hfill
+
(\ast)^T \Pi 
\left[
\begin{array}{ccc}
 1 & 0 & 0 \\
 \bin & \Win & 0 \\
 0 & 0 & I_n
\end{array}
\right]\preceq 0.  
\end{array}$}
\label{eq:main2}
\end{equation}
\label{th:main2}
\end{theorem}
\begin{proofof}{\rth{th:main2}}
For an input-output pair $(w,z)$ of the FNN $G$, 
let us define $q=\Win w+\bin$ and $p=\Phi(q)$.  
Then we have $z=\Wout p$.  
Multiplying $[\ 1\ w^T\ p^T\ ]^T$ from the right and its transpose from the left to 
\rec{eq:main2}, we obtain
\[
\scalebox{0.95}{$
\begin{array}{@{}l}
 -\Lsq+|z-z_0|_2^2+\tau(\varepsilon^2-|w-w_0|_2^2)
+ \left[
 \begin{array}{c}
 1 \\
 q \\
 p \\
 \end{array}
 \right]^T
 \Pi
 \left[
 \begin{array}{c}
 1 \\
 q \\
 p \\
 \end{array}
 \right]\le 0.  
\end{array}$}
\]
Since $\Pi\in\bPi^\star$, we readily obtain
\[
 -\Lsq+|z-z_0|_2^2+\tau(\varepsilon^2-|w-w_0|_2^2)\le 0.  
\]
Since $\tau\ge 0$, the above inequality
 implies that for any $w\in\clB_{\varepsilon}(w_0)$ 
we have $ -\Lsq+|z-z_0|_2^2\le 0$.  
This clearly shows that $L_{w_0,\varepsilon}\le \sqrt{\Lsq}$ holds.  
\end{proofof}

From \rth{th:main2}, we see that the upper bound computation 
of the Lipschitz constant for the FNN given by \rec{eq:FNN} can be formulated as
\begin{equation}
 \inf_{\Lsq,\tau,\Pi\in\bPi^\star}\ \Lsq\quad \mathrm{subject\ to}\ \rec{eq:main2}.  
\label{eq:primal}
\end{equation}

We note that \rth{th:main2} is a special case result of \cite{Fazlyab_IEEE2022}
where more general class of reliability verification problems are considered.  
Our novel contributions in this paper include:  
(i) providing novel copositive multipliers capturing accurately the 
input-output behavior of ReLUs;  
(ii) deriving an exactness verification test of the computed upper bounds by taking the dual of the problem \rec{eq:primal}; 
(iii) constructing reduced order models enabling us to deal with practical FNNs with 
hundreds of ReLUs.  

\section{Novel Multipliers for ReLUs and Comparison with Existing Ones}

In this section, we consider the concrete set of multipliers 
that satisfy the condition in \rec{eq:IQC0}.  
To this end, we define the positive semidefinite cone $\PSD^{n}\subset \bbS^n$, 
the copositive cone $\COP^{n}\subset \bbS^n$, and the nonnegative cone
$\NN^{n}\subset \bbS^n$ as follows:
\[
 \begin{array}{@{}l}
  \PSD^n:=\{P\in\bbS^n:\ x^TPx\geq 0\ \forall x\in\bbR^n\},\\
  \COP^n:=\{P\in\bbS^n:\ x^TPx\geq 0\ \forall x\in\bbR_{+}^n\},\\
  \NN^n:= \{P\in\bbS^n:\ P\geq 0\}.  
 \end{array}
\]
We can readily see that $\PSD^n\subset \PSD^n+\NN^n\subset \COP^n$.  
The copositive programming problem (COP) is a convex optimization problem in which 
we minimize a linear
objective function over the 
linear matrix inequality (LMI) constraints on the copositive cone   
\cite{Dur_2010}.  
As mentioned in \cite{Dur_2010}, 
the problem to determine whether a given symmetric matrix is
copositive or not is a co-NP complete problem.  
Therefore, it is hard to solve COP numerically in general.  
However, since the problem to determine whether a given matrix is in $\PSD+\NN$
can readily be reduced to an semidefinite programming problem (SDP), 
we can numerically solve the convex optimization problems on $\PSD+\NN$.  
Moreover, when $n\leq 4$, it is known that $\COP_n=\PSD_n+\NN_n$ 
hence those COPs with $n\leq 4$ can be reduced to SDPs.  

\subsection{Novel Copositive Multipliers for ReLU}

By focusing on \rec{eq:ReLU_alg}, 
we can the first main result of this paper.  
\begin{theorem}
Let us define $\bPi_\COP,\bPi_\NN\subset\bbS^{2n+1}$ by
\begin{equation}
\scalebox{0.95}{$
\begin{array}{@{}l}
\bPi_\COP:=
\left\{\Pi\in\bbS:\ 
\Pi= (\ast)^T  
\left(
Q+\clJ(J)
\right)E,\right.\\ 
\left.\hspace*{30mm}
J\in\bbD^n,\ Q\in\COP^{2n+1}
\right\},\\
\bPi_\NN:=
\left\{\Pi\in\bbS:\ 
\Pi= (\ast)^T  
\left(
Q+\clJ(J)
\right)E,\right.\\ 
\left.\hspace*{30mm}
 J\in\bbD^n,\ Q\in\NN^{2n+1}
\right\},\\
E:=\left[
\begin{array}{ccc}
 1 & 0_{1,n} & 0_{1,n}\\
 0_{n,1} & -I_n & I_n \\
 0_{n,1} & 0_n & I_n \\
\end{array}
\right],\ 
\clJ(J):=
\begin{bmatrix}
0 & 0_{1,n} & 0_{1,n} \\
\ast & 0_{n,n} & J \\
\ast & \ast & 0_{n,n} \\
\end{bmatrix}.  
\end{array}$}
\label{eq:PiCOP}
\end{equation}
Then we have $\bPi_\NN\subset\bPi_\COP\subset \bPi^\star$.  
\label{th:PiCOP}
\end{theorem}
\begin{proofof}{\rth{th:PiCOP}}
We see $\bPi_\NN\subset\bPi_\COP$ holds since $\NN\subset\COP$.  
Therefore it suffices to prove $\bPi_\COP\subset\bPi^\star$.  
For $\Pi\in\bPi_\COP$ and $p,q\in\bbR^n$ such that $p=\Phi(q)$, we have
\[
 \left[
\arraycolsep=0.5mm
 \begin{array}{c}
 1 \\
 q \\
 p \\
 \end{array}
 \right]^T\!\!\!
 \Pi
 \left[
\arraycolsep=0.5mm
 \begin{array}{c}
 1 \\
 q \\
 p \\
 \end{array}
 \right]=
 \left[
\arraycolsep=0.5mm
 \begin{array}{c}
 1 \\
 p-q \\
 p \\
 \end{array}
 \right]^T\!\!\!
 Q
 \left[
\arraycolsep=0.5mm
 \begin{array}{c}
 1 \\
 p-q \\
 p \\
 \end{array}
 \right]+2(p-q)^TJp.  
\]
Here, we see from \rec{eq:ReLU_alg}, $Q\in\COP^{2n+1}$, and $J\in\bbD^n$ that
\[
 \left[
\arraycolsep=0.5mm
 \begin{array}{c}
 1 \\
 p-q \\
 p \\
 \end{array}
 \right]^T\!\!\!
 Q
 \left[
\arraycolsep=0.5mm
 \begin{array}{c}
 1 \\
 p-q \\
 p \\
 \end{array}
 \right]\ge 0,\ (p-q)^TJp=0.  
\]
It follows that
\[
 \left[
\arraycolsep=0.5mm
 \begin{array}{c}
 1 \\
 q \\
 p \\
 \end{array}
 \right]^T\!\!\!
 \Pi
 \left[
\arraycolsep=0.5mm
 \begin{array}{c}
 1 \\
 q \\
 p \\
 \end{array}
 \right]\ge 0.  
\]
This clearly shows $\Pi\in\bPi^\star$.  
\end{proofof}

In the following, we review existing multipliers capturing the behavior of ReLUs.  
Then, we clarify the inclusion relationship among those existing ones and
the set of copositive multipliers $\bPi_\COP$ and its inner approximation $\bPi_\NN$.  

\subsection{O'Shea-Zames-Falb Multipliers}

We now review the arguments in \cite{Fetzler_IFAC2017} on 
O'Shea-Zames-Falb multipliers \cite{O'Shea_IEEE1967,Zames_SIAM1968}.  
The overview of O'Shea-Zames-Falb multipliers can be found at \cite{Carrasco_EJC2016}, 
where the contributions of O'Shea to the development of the multipliers are emphasized.  
Due to this reason, we call the multipliers as O'Shea-Zames-Falb multipliers
as opposed to the well celebrated moniker Zames-Falb multipliers.    
This modified moniker has already been employed, e.g., in \cite{Scherer_Arxiv2022}.  

We first introduce the following definition.  
\begin{definition}
Let $\mu\le 0 \le \nu$.  
Then the nonlinearity $\phi:\ \bbR\to\bbR$ is slope-restricted, 
in short $\phi \in \slope(\mu,\nu)$, if $\phi(0)=0$ and
\[
  \mu\leq \dfrac{\phi(x)-\phi(y)}{x-y}
 \le
  \sup_{x\neq y}\dfrac{\phi(x)-\phi(y)}{x-y}
  <\nu
\]
 for all $x,y\in\bbR$, $x\neq y$.
\end{definition}

The main result of \cite{Fetzler_IFAC2017}
on the static O'Shea-Zames-Falb multipliers for slope-restricted nonlinearities can be
summarized by the next lemma.
\begin{lemma}\cite{Fetzler_IFAC2017}
For a given nonlinearity $\phi\in \slope(\mu,\nu)$ with 
$\mu\le 0\le \nu$, 
let us define $\Phi:\ \bbR^m\to \bbR^m$
by $\Phi:=\diag(\phi,\cdots,\phi)$.  
Assume $M\in\DHD^m$.  Then we have
\[
\scalebox{0.93}{$
\arraycolsep=0.5mm
\begin{array}{@{}l}
 (\ast)^T
\left[
 \begin{array}{cc}
  0 & M^T\\
  M & 0 
 \end{array}
\right]
\left(
 \left[
 \begin{array}{cc}
  \nu I_m & -I_m\\
  -\mu I_m & I_m\\
 \end{array}
 \right]
 \left[
 \begin{array}{c}
  x\\
  \Phi(x)
 \end{array}
 \right]\right)\ge 0\ \forall x \in\bbR^m.    
\end{array}$}
\]
\label{le:sl}
\end{lemma}

From this key lemma and the fact that the ReLU $\phi:\ \bbR\to\bbR$
satisfies $\phi\in\slope(0,1)$, 
we can obtain the next result on 
the O'Shea-Zames-Falb multipliers for the ReLU given by \rec{eq:ReLU}.  
\begin{proposition}
Let us define
\[
\scalebox{0.83}{$
\begin{array}{@{}l}
\arraycolsep=0.5mm
 \hPiZF:=
 \left\{
  \hPi\in\bbS:\ 
  \hPi=
 (\ast)^T
\left[
 \begin{array}{cc}
  0 & M\\
  M^T & 0 
 \end{array}
\right]
 \left[
 \begin{array}{cc}
  I_n & -I_n\\
  0 & I_n\\
 \end{array}
 \right],\ M\in \DHD^n \right\},
\end{array}$}
\]
\begin{equation}
\scalebox{0.83}{$
\begin{array}{@{}l}
 \PiZF:=
 \left\{
  \Pi\in\bbS:\ 
  \Pi=\diag(0,\hPi),\ \hPi\in\hPiZF \right\}.  
\end{array}$}
\label{eq:MOZF}
\end{equation}
Then we have $\PiZF\subset\bPi^\star$.  
\label{pro:OZF}
\end{proposition}
%

\subsection{Multipliers by Fazlyab et al. \cite{Fazlyab_IEEE2022}}

We next review the multipliers introduced by Fazlyab et al. in \cite{Fazlyab_IEEE2022}.  
To this end, let us define
\[
\scalebox{0.85}{$
\begin{array}{@{}l}
\arraycolsep=0.5mm
\displaystyle
\bbT^n:=
\left\{
T\in\bbS^n:\ T=\sum_{1\le i < j \le n}\lambda_{ij}(e_i-e_j)(e_i-e_j)^T,\ \lambda_{ij}\ge 0
\right\}.    
\end{array}$}
\]
Then, the result in Lemma 3 of \cite{Fazlyab_IEEE2022} can be summarized 
as follows. 
\begin{proposition}\cite{Fazlyab_IEEE2022}
Let us define $\PiFaz\in\bbS^{2n+1}$ by
\begin{equation}
\begin{array}{@{}ll}
 \PiFaz:=&\left\{
	 \Pi\in\bbS:\ \Pi=
	 \begin{bmatrix}
	  0 & -\nu & \nu+\eta \\
	  \ast & 0 & \Lambda + T\\
	  \ast & \ast & -2(\Lambda+T)
	 \end{bmatrix},\right.\\ 
         &\nu^T,\eta^T\in\bbR^n_+,\ \Lambda\in\bbD^n,\ T\in\bbT^n \Biggr \}.  
\end{array}
 \label{eq:MFaz} 
\end{equation}
Then we have $\PiFaz\subset\bPi^\star$.  
\end{proposition}
%

\subsection{Comparison among Multipliers}

We are now ready to clarify the relationship among the set of multipliers 
$\bPi_\COP$ and $\bPi_\NN$ given by \rec{eq:PiCOP}, 
$\PiZF$ given by \rec{eq:MOZF}, and 
$\PiFaz$ given by \rec{eq:MFaz}.  
The next theorem is one of the main results of the present paper.  
\begin{theorem}
For the set of multipliers 
$\bPi_\COP$ and $\bPi_\NN$ given by \rec{eq:PiCOP}, 
$\PiZF$ given by \rec{eq:MOZF}, and 
$\PiFaz$ given by \rec{eq:MFaz}, we have
\begin{equation}
\PiZF\subset \bPi_\NN \subset \bPi_\COP,
\label{eq:inclOZF}
\end{equation}
\begin{equation}
\PiFaz\subset \bPi_\NN \subset \bPi_\COP.  
\label{eq:inclFaz}
\end{equation}
\label{th:incl}
\end{theorem}
\begin{proofof}{\rth{th:incl}}
For the proof, we note that $\bPi_\NN$ given by \rec{eq:PiCOP} and
$\PiZF$ given by \rec{eq:MOZF} can be rewritten, equivalently, as follows:
\[
\scalebox{0.70}{$
 \begin{array}{@{}l}
  \bPi_\NN=
\left\{\Pi\in\bbS:\ 
\Pi= \begin{bmatrix}
   Q_{11}  & -Q_{12}  & Q_{12}+Q_{13} \\
   \ast  & Q_{22}  & -Q_{22}-Q_{23}-J \\
   \ast & \ast  & Q_{22}+ Q_{33}+Q_{23}+Q_{23}^T+2J   
     \end{bmatrix},\ 
 J\in\bbD^n,\ \right.
\end{array}$}
\]
\begin{equation}
\scalebox{0.70}{$
 \begin{array}{@{}l}
\left.
 Q=\begin{bmatrix}
   Q_{11}  & Q_{12}  & Q_{13} \\
   \ast  & Q_{22}  & Q_{23} \\
   \ast & \ast  & Q_{33} 
   \end{bmatrix}\in\NN^{2n+1},\ Q_{11}\in\bbR_+,\ Q_{22},Q_{33}\in\NN^n
\right\},
\end{array}$}
\label{eq:COPexp}
\end{equation}
\begin{equation}
\scalebox{0.70}{$
 \begin{array}{@{}l}
  \PiZF:=
 \left\{
  \Pi\in\bbS:\ 
\Pi= \begin{bmatrix}
   0  & 0_{1,n}  & 0_{1,n} \\
   \ast  & 0_{n,n}  & M \\
   \ast & \ast  & -M-M^T
     \end{bmatrix},\ M\in \DHD^n \right\}.  
 \end{array}$}
\label{eq:ZFexp}
\end{equation}
We first prove \rec{eq:inclOZF}.  To this end, it suffices to prove $\PiZF\subset \bPi_\NN$.  
For $\Pi\in\PiZF$ of the form \rec{eq:ZFexp}, 
let us choose $Q$ in \rec{eq:COPexp} as
$Q_{11}=0$, $Q_{12}=Q_{13}=0$, $Q_{22}=Q_{33}=0$, and select 
$Q_{23}$ and $J$ such that $Q_{23}+J=-M$ that is always feasible.  
Then we see that $\Pi\in\bPi_\NN$ and hence $\PiZF\subset \bPi_\NN$ holds.  
We next prove \rec{eq:inclFaz}.  To this end, it suffices to prove $\PiFaz\subset \bPi_\NN$.  
For $\Pi\in\PiFaz$ of the form \rec{eq:MFaz}, 
let us choose $Q$ in \rec{eq:COPexp} as
$Q_{11}=0$, $Q_{12}=\nu$, $Q_{13}=\eta$, $Q_{22}=Q_{33}=0$, and select 
$Q_{23}$ and $J$ such that $Q_{23}+J=-(\Lambda+T)$ that is always feasible.  
Then we see that $\Pi\in\bPi_\NN$ and hence $\PiFaz\subset \bPi_\NN$ holds.  
\end{proofof}

On the basis of this result, we now describe the SDP
for the upper bound computation of 
the local Lipschitz constant of the FNN given by \rec{eq:FNN} and \rec{eq:ReLU}.  
Since the inclusion relationships in \rth{th:incl} hold, 
and since $\bPi_\COP$ is numerically intractable, 
we focus on the following SDP:   
\begin{equation}
\gamma_\primal^2:= \inf_{\Lsq,\tau,\Pi\in\bPi_\NN}\  \Lsq\ \mathrm{subject\ to}\ \rec{eq:main2}.  
\label{eq:primal_SDP}
\end{equation}
If this SDP is feasible, we can conclude that $L_{w_0,\varepsilon}\le \gamma_\primal$.  

\begin{remark}
The usefulness of the copositive multipliers or its inner approximation
has been already observed in the stability analysis of 
recurrent neural networks (RNNs) with activation functions being ReLUs, 
see \cite{Ebihara_EJC2021,Ebihara_CDC2021}.  
\end{remark}

\section{Dual SDP and Exactness Verification}

By solving the (primal) SDP \rec{eq:primal_SDP}, 
we can obtain an upper bound of the local Lipschitz constant.  
However, if we merely rely on such upper bound computation, 
we cannot draw any definite conclusion on how far the upper bound is close to the exact one.  
To address this issue, it is well known in the field of robust control that
considering the dual problem is useful \cite{Scherer_SIAM2005,Ebihara_IEEE2009}.  
By following the Lagrange duality theory for SDPs \cite{Scherer_EJC2006},  
the dual SDP of \rec{eq:primal_SDP} can be obtained as follows:  
\begin{equation}
\scalebox{0.95}{$
\begin{array}{@{}l}\displaystyle
\gamma_\dual^2:=\sup_{H\in\bbS_+^{1+m+n}} 
\trace\left(
(\ast)^T H 
\left[
 \begin{array}{c}
z_0^T \\
0_{m,l}  \\
-\Wout^T
 \end{array}
\right]\right)\ \mathrm{subject\ to}\hspace*{-20mm}\\
H_{11}=1,\\ 
\trace\left(
\begin{bmatrix}
\varepsilon^2-w_0^Tw_0 & w_0^T \\
\ast & - I_m 
\end{bmatrix}
\begin{bmatrix}
H_{11} & H_{12} \\
\ast & H_{22} 
\end{bmatrix}
\right)\ge 0,\\
\diag(-\bin H_{13}-\Win H_{23}+H_{33})=0,\\
(\ast)^TH
\left[
\begin{array}{ccc}
 1 & 0_{1,m} & 0_{1,n}\\
 -\bin & -\Win & I_n \\
 0_{1,n} & 0_{n,m} & I_n \\
\end{array}
\right]^T \ge 0,\\
 H =:
\begin{bmatrix}
 H_{11} & H_{12} & H_{13}\\
 \ast & H_{22} & H_{23}\\
 \ast & \ast & H_{33}
\end{bmatrix},\ 
H_{22}\in \bbS_+^m,\ 
H_{33}\in \bbS_+^n
\end{array}$}
\label{eq:dual_SDP}
\end{equation}

We can prove that the primal SDP \rec{eq:primal_SDP} is strongly feasible.  
Therefore, the dual SDP \rec{eq:dual_SDP} has an optimal solution, 
and there is no duality gap between \rec{eq:primal_SDP} and \rec{eq:dual_SDP} and thus
$\gamma_\primal=\gamma_\dual$ holds \cite{Klerk_2002}.  
By focusing on the dual SDP \rec{eq:dual_SDP}, we can obtain the next results
on the exactness verification of the computed upper bounds.  
\begin{theorem}
Suppose $\rank(H)=1$ holds for the optimal solution $H\in\bbS_+^{1+m+n}$
of the dual SDP \rec{eq:dual_SDP}.  
Then we have $L_{w_0,\varepsilon} =\gamma_\dual$.  
Moreover, the full-rank factorization of the optimal solution $H\in\bbS_+^{1+m+n}$
is given by 
\begin{equation}
 H=
 \begin{bmatrix}
    1 \\ h_2 \\ h_3 
 \end{bmatrix}
 \begin{bmatrix}
    1 \\ h_2 \\ h_3 
   \end{bmatrix}^T,\ h_2\in\bbR^m,\ h_3\in\bbR^n.   
\label{eq:FF}
\end{equation}
In addition, if we define $w^\star:=h_2$, then $w^\star\in\clB_{\varepsilon}(w_0)$ holds, and 
$w^\star\in\clB_{\varepsilon}(w_0)$ is one of the worst-case inputs
satisfying $L_{w_0,\varepsilon} =|G(w^\star)-G(w_0)|_2$.  
\label{th:EV}
\end{theorem}
\begin{proofof}{\rth{th:EV}}
Suppose $\rank(H)=1$ holds for the optimal solution $H\in\bbS_+^{1+m+n}$.  
Then, since $H_{11}=1$, we readily see that the full-rank factorization of $H$ is given of the form of \rec{eq:FF}.  
With this rank-one factorization and the form of the dual SDP \rec{eq:dual_SDP}, we have
\begin{equation}
 \gamma_\dual^2=|z_0-\Wout h_3|_2^2,
\label{eq:dual1}
\end{equation}
\begin{equation}
 \varepsilon^2\ge |w_0-h_2|_2^2,
\label{eq:dual2}
\end{equation}
\begin{equation}
 \diag\left((h_3-(\Win h_2+\bin))h_3^T\right)= 0,
\label{eq:dual3}
\end{equation}
\begin{equation}
 \begin{bmatrix}
 1\\ h_3-(\Win h_2+\bin)\\ h_3 
 \end{bmatrix}
 \begin{bmatrix}
 1\\ h_3-(\Win h_2+\bin)\\ h_3 
 \end{bmatrix}^T\ge 0.  
\label{eq:dual4}
\end{equation}
Here we can rewrite \rec{eq:dual3} equivalently as 
\begin{equation}
 (h_3-(\Win h_2+\bin))\otimes h_3= 0.  
\label{eq:dual7}
\end{equation}
From \rec{eq:dual4}, \rec{eq:dual7}, and \rec{eq:ReLU_alg}, we obtain 
\begin{equation}
 h_3=\Phi(\bin+\Win h_2).  
\end{equation}
From this equation together with \rec{eq:dual1} and \rec{eq:dual2}, we see that for 
$w^\star:=h_2\in\clB_{\varepsilon}(w_0)$ we have 
$\gamma_\dual=|G(w^\star)-G(w_0)|_2$.  
From this fact and $\gamma_\dual\ge L_{w_0,\varepsilon}$, it follows that 
$\gamma_\dual= L_{w_0,\varepsilon}$.  
In addition, we can conclude that 
$w^\star\in\clB_{\varepsilon}(w_0)$ is one of the worst-case inputs achieving
$L_{w_0,\varepsilon}=|G(w^\star)-G(w_0)|_2$. 
\end{proofof}
\begin{remark}
If we employ the set of multipliers such as 
$\PiZF$ and $\PiFaz$ in the primal SDP \rec{eq:primal_SDP}, 
we are led to a dual SDP that is of course different from \rec{eq:dual_SDP}.   
In such a case, we cannot draw the exactness verification condition shown in \rth{th:EV}.  
This in part verifies the soundness of employing the new set of multipliers $\bPi_\NN$.  
\end{remark}
\begin{remark}
In another problem setting for certifying robustness of FNNs, 
the exactness verification of SDP relaxations is discussed in \cite{Zhang_arxiv}.  
\end{remark}
%
\section{Exact Model Reduction around Target Input}
\label{sec:reduction}

When solving the SDP \rec{eq:primal_SDP}, 
we have to deal with the multiplier variable 
$\Pi\in\bPi^\star\subset\bbS^{2n+1}$.  
Here, $n$ stands for the number of ReLUs and this would be
at least a few hundred in practical FNNs.  
Therefore, in practical problem settings the SDP \rec{eq:primal_SDP}
becomes intractable.   
To address this issue, in this section, 
we propose a method for the reduction of the number of ReLUs
while maintaining the input-output behavior of the original FNN, 
by making full use of the property of ReLUs and the local Lipschitz constant computation 
problem.  

To be precise, let us consider the behavior of 
\[
\begin{array}{@{}lcl}
 z&=&\Wout\Phi(\Win w+\bin),\ \Win\in\bbR^{n\times m},\ \Wout\in\bbR^{l\times n} 
\end{array}
\]
for $w\in  \clB_{\varepsilon}(w_0)$. To this end, let us define 
$q=\Win w+\bin$ and $q_0=\Win w_0+\bin$.  Then we see
\[
 \min_{w\in \clB_{\varepsilon}(w_0)} q_i=q_{0,i}-\varepsilon |\Win_{,i}|_2,\ 
 \max_{w\in \clB_{\varepsilon}(w_0)} q_i=q_{0,i}+\varepsilon |\Win_{,i}|_2
\]
where $\Win_{,i}\in\bbR^{1\times m}$ stands for 
the $i$-th row of $\Win\in\bbR^{n\times m}$.  
It follows that 
\begin{equation}
 \left\{
\begin{array}{lcl}
 q_{0,i}\ge \varepsilon |\Win_{,i}|_2 & \Ra & q_i \ge 0\ \forall w\in  \clB_{\varepsilon}(w_0), \\
 q_{0,i}\le -\varepsilon |\Win_{,i}|_2  & \Ra &  q_i \le 0\ \forall w\in  \clB_{\varepsilon}(w_0).    
\end{array}
\right.  
\label{eq:ReLU_red}
\end{equation}
With this fact in mind, let us define
\begin{equation}
 \begin{array}{@{}l}
  \clZ_n:=\{1,2,\cdots,n\},\\
  \clN_+:=\{i\in \clZ_n:\ q_{0,i}\ge \varepsilon |\Win_{,i}|_2 \},\\
  \clN_0:=\{i\in \clZ_n:\ q_{0,i}\le -\varepsilon |\Win_{,i}|_2 \},\\
  \clN_r:=\clZ_n\setminus\{\clN_+\cup\clN_0\}.    
 \end{array}
\label{eq:red_index}
\end{equation}
Then, 
$q_i$ with $i\in\clN_+$ is never rectified for all $w\in\clB_{\varepsilon}(w_0)$, 
$q_i$ with $i\in\clN_0$ is rectified for all $w\in\clB_{\varepsilon}(w_0)$, and
for $q_i$ with $i\in\clN_r$ we cannot say anything definitely.  
This motivates us to define
\begin{equation}
 \begin{array}{@{}l}\displaystyle
 E_+:=\bigoplus_{i\in\clN_+} e_i,\  E_r:=\bigoplus_{i\in\clN_r} e_i,\\ 
 \tWin:=E_+\Win,\ \tbin:=E_+\bin,  \tWout:=\Wout E_+^T, \\
 \hWin:=E_r\Win,\ \hbin:=E_r\bin,  \hWout:=\Wout E_r^T.  
 \end{array}
\label{eq:red_coef}
\end{equation}
Then, we can obtain the next reduced order model 
with the number of ReLUs being $n_r:=|\clN_r|$.  
\begin{equation}
G_r: z=\tWout(\tWin w+\tbin)+\hWout\Phi(\hWin w+\hbin).  
\label{eq:FNN_red}
\end{equation}
By the construction procedure of this reduced order model, the next results readily hold.  
\begin{theorem}
For the FNN $G$ given by \rec{eq:FNN} and \rec{eq:ReLU}, 
let us consider the reduced order model $G_r$ constructed by
\rec{eq:FNN_red}, \rec{eq:red_index}, and \rec{eq:red_coef}.  
Then, we have
\begin{equation}
G(w)=G_r(w) \ \forall w\in\clB_\varepsilon(w_0).  
\end{equation}
\end{theorem}
\begin{remark}
From the observation \rec{eq:ReLU_red}, 
we have constructed explicitly an exact reduced order model $G_r$
whose input-output behavior is identical to the original FNN $G$
for inputs $w\in\clB_\varepsilon(w_0)$.  
Similar observations to \rec{eq:ReLU_red}
are used in \cite{Fazlyab_IEEE2022} for 
the purpose of tightening multiplier relaxations.  
\end{remark}

By using the reduced order model $G_r$, 
we can derive the primal and dual SDPs 
for the upper bound computation of the local Lipschitz constant
$L_{w_0,\varepsilon}$ of the original FNN $G$.    
These are given as follows:  

\noindent Primal SDP:  
\begin{equation}
\scalebox{0.83}{$
 \begin{array}{@{}l}\displaystyle
\gamma_\rprimal^2:=\inf_{\Lsq, \tau, \Pi\in\bPi_\NN}\ \Lsq\ 
\mathrm{subject\ to}\\
\arraycolsep=1mm
(*)^T  \!\!
\left[
\begin{array}{ccc} 
-\Lsq+\tau\epsilon^2 & 0 & 0\\
\ast & -\tau I_m & 0\\
\ast & \ast & I_l
\end{array}
\right]\!\!\!
\left[
\begin{array}{ccc} 
1 & 0 & 0\\ 
-w_0 & I_m & 0\\ 
\tWout \tbin-z_0 & \tWout \tWin & \hWout 
\end{array}\right] \hspace*{-20mm}\\
\hfill +(\ast)^T \Pi
\begin{bmatrix}
1 & 0 & 0 \\
\hbin & \hWin & 0 \\
0 & 0 & I_{n_r}
\end{bmatrix}
\preceq 0.  
 \end{array}$}
\label{eq:primal_SDP_red}
\end{equation}

\noindent Dual SDP:  
\begin{equation}
\scalebox{0.83}{$
\begin{array}{@{}l}\displaystyle 
\gamma_\rdual^2:=\!\!\!\!\sup_{H\in\bbS_+^{1+m+{n_r}}}\!\!\!\!
\trace\left(
(\ast)^T H 
\left[
 \begin{array}{c}
(\tWout\tbin-z_0)^T \\
\tWin^T\tWout^T  \\
\hWout^T
 \end{array}
\right]\right)\ \mathrm{subject\ to}\hspace*{-20mm}\\
H_{11}=1,\\ 
\trace\left(
\begin{bmatrix}
\varepsilon^2-w_0^Tw_0 & w_0^T \\
\ast & - I_m 
\end{bmatrix}
\begin{bmatrix}
H_{11} & H_{12} \\
\ast & H_{22} 
\end{bmatrix}
\right)\ge 0,\\
\diag(-\hbin H_{13}-\hWin H_{23}+H_{33})=0,\\
(\ast)^TH
\left[
\begin{array}{ccc}
 1 & 0_{1,m} & 0_{1,{n_r}}\\
 -\hbin & -\hWin & I_{n_r} \\
 0_{1,{n_r}} & 0_{{n_r},m} & I_{n_r} \\
\end{array}
\right]^T \ge 0,\\
H = 
\begin{bmatrix}
 H_{11} & H_{12} & H_{13}\\
 \ast & H_{22} & H_{23}\\
 \ast & \ast & H_{33}
\end{bmatrix},\  
H_{22}\in \bbS_+^m,\ 
H_{33}\in \bbS_+^{n_r}.  
\end{array}$}
\label{eq:dual_SDP_red}
\end{equation}

In the primal SDP \rec{eq:primal_SDP_red}, we emphasize that 
the size of the multiplier variable $\Pi$ has been reduced to $2n_r+1$.  
Similarly, in the dual SDP \rec{eq:dual_SDP_red}, 
the size of the dual variable $H$ has been reduced to $1+m+n_r$.  

Again, we can prove that the primal SDP \rec{eq:primal_SDP_red} is strongly feasible.  
Therefore, the dual SDP \rec{eq:dual_SDP_red} has an optimal solution, 
and there is no duality gap between \rec{eq:primal_SDP_red} and 
\rec{eq:dual_SDP_red}.  Namely, 
$\gamma_\rprimal=\gamma_\rdual\ge L_{w_0,\varepsilon}$ holds.  
In addition, as for the exactness verification, the next results follow.  
\begin{theorem}
Suppose $\rank(H)=1$ holds for the optimal solution $H\in\bbS_+^{1+m+n_r}$
of the dual SDP \rec{eq:dual_SDP_red}.  
Then we have $L_{w_0,\varepsilon} =\gamma_\rdual$.  
Moreover, the full-rank factorization of the optimal solution $H\in\bbS_+^{1+m+n_r}$
is given by 
\begin{equation}
 H=
 \begin{bmatrix}
    1 \\ h_2 \\ h_3 
 \end{bmatrix}
 \begin{bmatrix}
    1 \\ h_2 \\ h_3 
   \end{bmatrix}^T,\ h_2\in\bbR^m,\ h_3\in\bbR^{n_r}.   
\label{eq:FFred}
\end{equation}
In addition, if we define $w^\star:=h_2$, then $w^\star\in\clB_{\varepsilon}(w_0)$ holds, and 
$w^\star\in\clB_{\varepsilon}(w_0)$ is one of the worst-case inputs
satisfying $L_{w_0,\varepsilon} =|G(w^\star)-G(w_0)|_2$.  
\label{th:EVred}
\end{theorem}
\begin{remark}
There is no magnitude relationship between 
$\gamma_\primal=\gamma_\dual$ and
$\gamma_\rprimal=\gamma_\rdual$.  
Still, in practical problem settings with hundreds of ReLUs,
the computation of $\gamma_\primal=\gamma_\dual$ becomes intractable, 
whereas the computation of $\gamma_\rprimal=\gamma_\rdual$ remains to be tractable.  
We illustrate this point, as well as the usefulness of the exactness verification method, 
by numerical examples in the next section.  
\end{remark}
%

\section{Numerical Examples}

In this section, we illustrate the usefulness of the model reduction
and exactness verification methods by numerical examples.  
When solving SDPs, 
we used MATLAB 2023a and 
the SDP solver MOSEK \cite{MOSEK} together with the 
parser YALMIP \cite{Lofberg_IEEECACSD2004} 
on a computer with CPU 12th Gen Intel(R) Core(TM) i9-12900 2.40 GHz.

\subsection{Academic Toy Examples}

Let us consider the case where
$n=6$, $m=3$, $l=3$ in \rec{eq:FNN} and
\[
\begin{array}{@{}l}
\Win=
\left[
\begin{array}{rrr}
   -0.62 &  -0.28 &   0.47\\
    0.88 &   0.18 &   0.48\\
    0.37 &  -0.12 &   0.40\\
    0.22 &   0.16 &   0.10\\
    0.31 &   0.90 &   0.49\\
    0.42 &   0.39 &  -0.56\\
 \end{array}
\right],\ 
\bin=
\left[
\begin{array}{r}
   -0.18\\
    0.71\\
    0.34\\
   -0.09\\
    0.22\\
   -0.20\\
 \end{array}
\right],\\
\Wout=
\left[
\begin{array}{rrrrrr}
    0.19 &   0.30 &   0.38 &   0.51 &  -0.79 &  -0.74\\
    0.35 &   0.12 &   0.07 &   0.39 &   0.42 &  -0.18\\
    0.00 &  -0.62 &  -0.14 &  -0.60 &   0.04 &   0.47\\
 \end{array}
\right].  
\end{array}
\]
where these are randomly generated.  
As for the target input, we randomly chose
\[
 w_0=
\left[
\begin{array}{rrr}   
0.52 &  -0.15  &  -0.07
 \end{array}
\right]^T,\ |w_0|_2=0.5457.    
\]
The corresponding output turns out to be
\begin{equation}
 z_0=
\left[
\begin{array}{rrr}
   0.3632 &  0.2584 &  -0.7510
 \end{array}
\right]^T.  
\label{eq:z0}
\end{equation}

By letting $\varepsilon=0.1$, we solved the primal SDP \rec{eq:primal_SDP}  
and obtained $\gamma_\primal=0.1088$.  
The CPU time was 0.053 [sec].  
In addition, by the command \texttt{dual} implemented in YALMIP, 
we extracted the dual optimal solution $H$, 
and its rank was numerically verified to be one.  
By following \rth{th:EV}, we constructed $w^\star\in\bbR^3$ as 
\[
w^\star=
\left[
\begin{array}{rrr}   
 0.5115 &  -0.0648 &  -0.1217
 \end{array}
\right]^T.  
\]
It turned out that
\[
 |w^\star-w_0|_2=0.1000,\ |G(w^\star)-G(w_0)|_2=0.1088.   
\]
Namely, the computed upper bound is exact and hence $L_{w_0,\varepsilon}=0.1088$ 
holds.  In addition, $w^\star\in\bbR^3$ is one of the worst-case inputs achieving
$L_{w_0,\varepsilon}=|G(w^\star)-G(w_0)|_2$.  

\subsection{Practical Real Examples}

In this section, we use the MNIST classifier (SDP-NN) described in 
\cite{Raghunathan_ICLR2018}, which is a single-layer ReLU-FNN.   
The inputs are $28\times 28$ pixel data of handwritten digits from $0$ to $9$.  
The magnitude of each pixel is normalized to $[0\ 1]$.  
The number of ReLUs used in this FNN is $500$.  
It follows that $m=784$, $l=10$, and $n=500$ in \rec{eq:FNN}.  
The classification rule is given by $C(w):=\argmax_{1\le i \le 10} G_i(w)-1$.  

\begin{figure}[t]
\begin{center}
\begin{minipage}{0.3\linewidth}
\centering
 \includegraphics[scale=0.16]{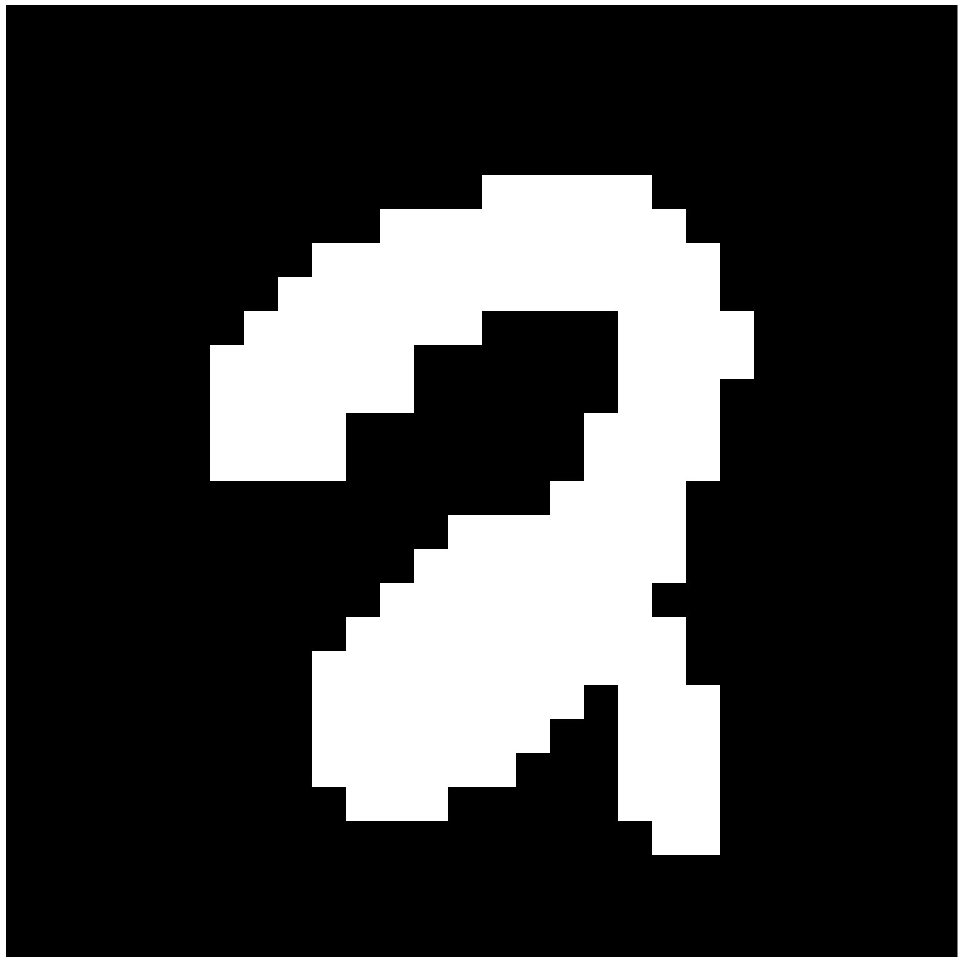}
\end{minipage}
\begin{minipage}{0.6\linewidth}
\quad $\Ra$\quad  FNN\ $\Ra$ 
$
 z_0 =
\begin{bmatrix*}[r]
  0.1506\\
   -0.1238\\
    \mathbf{0.7448}\\
   -0.1835\\
   -0.1983\\
   -0.2167\\
    0.0474\\
    0.3736\\
    0.5347\\
    0.3528\\
\end{bmatrix*}$  
\end{minipage}
 \caption{Image of the input $w_0$ and corresponding output $z_0=G(w_0)$.  }\vspace*{-5mm}
\label{fig:w0z0}
\end{center}
\end{figure}

As for the target input $w_0\in\bbR^m$, we chose the one whose image is shown
in \rfig{fig:w0z0} where $|w_0|_2=9.9652$.  
The corresponding output $z_0=G(w_0)$ is also shown in \rfig{fig:w0z0}.  
We can confirm that the input $w_0$ is correctly classified as $2$.  

\begin{figure}[t]
\centering
 \includegraphics[scale=0.62]{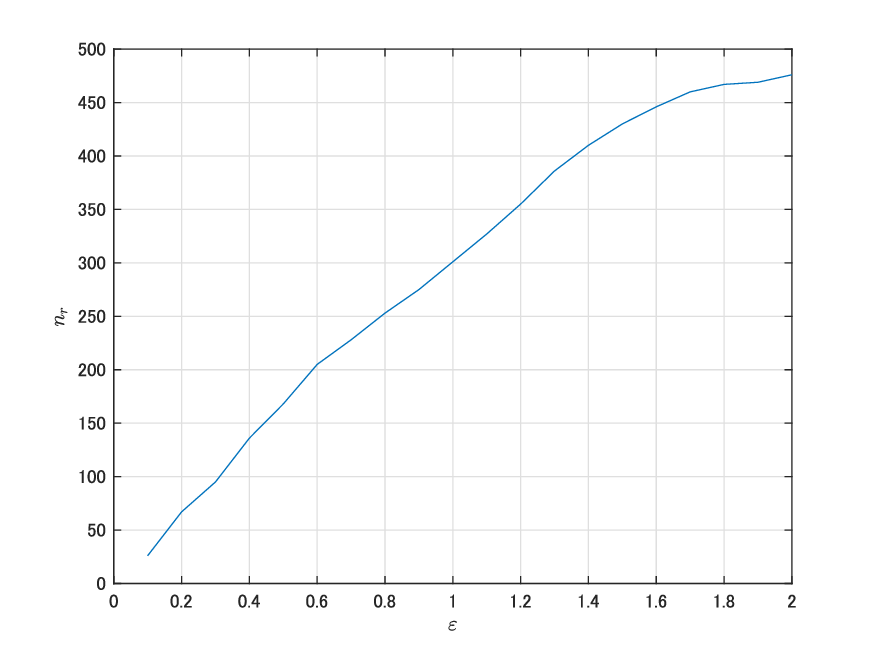}\vspace*{-3mm}
 \caption{Model reduction results.  }\vspace*{-5mm}
 \label{fig:reduction}
\end{figure}

If we directly work on this FNN and solve the primal SDP \rec{eq:primal_SDP}, 
we have to deal with the multiplier variable $\Pi$ of size $2n+1=1001$.  
However, this is numerically intractable.  
We therefore applied the exact model reduction method described 
in \rsec{sec:reduction}.  The results are shown in \rfig{fig:reduction}.  
As expected, we have achieved considerable reduction of the numbers of ReLUs
especially when $\varepsilon$ is small.  
In particular, when $\varepsilon=0.1$, the number of ReLUs in 
the reduced order model $G_r$ given by \rec{eq:FNN_red} was $n_r=26$.  

\begin{figure}[t]
\begin{center}
\begin{minipage}{0.3\linewidth}
\centering
 \includegraphics[scale=0.29]{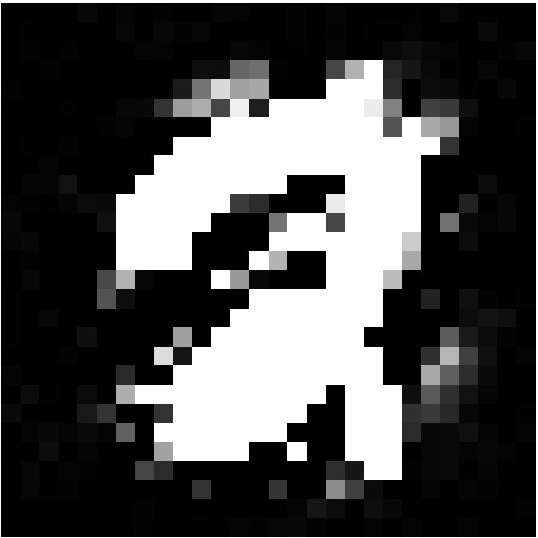}
\end{minipage}
\begin{minipage}{0.6\linewidth}
\quad $\Ra$\quad  FNN\ $\Ra$ 
$
 z^\star =
\begin{bmatrix*}[r]
    0.1017\\
   -0.1703\\
\mathbf{0.6991}\\
   -0.2244\\
   -0.2384\\
   -0.2613\\
    0.0026\\
    0.3275\\
    0.4931\\
    0.3051\\
\end{bmatrix*}$  
\end{minipage}
 \caption{Image of the input $w^\star$ and corresponding output $z^\star=G(w^\star)$.  }\vspace*{-5mm}
\label{fig:wwczwc}
\end{center}
\end{figure}

By letting $\varepsilon=0.1$, we solved the primal SDP \rec{eq:primal_SDP_red}
relying on the reduced order model and 
obtained $\gamma_\rprimal=0.1416$.  
The CPU time was 2289 [sec].  
In addition, by the command \texttt{dual}, 
we extracted the dual optimal solution $H$, 
and its rank was numerically verified to be one.  
From \rth{th:EVred}, we then constructed $w^\star\in\bbR^{784}$ 
whose image is shown in \rfig{fig:wwczwc}.  
The corresponding output $z^\star=G(w^\star)$ is also shown in \rfig{fig:wwczwc}.  
It turned out that
\[
 |w^\star-w_0|_2=0.1000,\ |G(w^\star)-G(w_0)|_2=0.1416.   
\]
Namely, the computed upper bound is exact and hence $L_{w_0,\varepsilon}=0.1416$ 
holds.  In addition, $w^\star\in\bbR^{784}$ is one of the worst-case inputs achieving
$L_{w_0,\varepsilon}=|G(w^\star)-G(w_0)|_2$.  
Here we emphasize that the $w_0$ is a genuine worst-case input for the original FNN $G$, 
even though it is constructed from the reduced order model $G_r$.  

We see from \rfig{fig:wwczwc} that even the worst-case input $w^\star$
is correctly classified as $2$.  This assertion can be strengthened 
by the robustness test \rec{eq:robl2cond} in \rpr{pr:rob}.  
Namely, regarding $z_0=G(w_0)$ shown in \rfig{fig:w0z0}, we see that 
\[
 L_{w_0,\varepsilon}= 0.1416 < 0.1485 = \frac{1}{\sqrt{2}}(0.7448-0.5347).  
\]
Therefore we are led to the definite conclusion that 
there is no adversarial input that leads to false-classification 
within $w\in\clB_{w_0,\varepsilon}$ for $\varepsilon=0.1$.  

%

\section{Conclusion and Future Works}

In this paper, we considered the computation problem of 
the local Lipschitz constants for single-layer FNNs.  
By following standard procedure using multipliers, 
we reduced the upper bound computation problem into an SDP.  
Our novel contributions in this paper include:  
(i) providing novel copositive multipliers capturing accurately the 
input-output behavior of ReLUs;  
(ii) deriving an exactness verification test of the computed upper bounds by taking the dual of the SDP; 
(iii) constructing reduced order models enabling us to deal with practical FNNs 
with hundreds of ReLUs.  
We finally illustrated the usefulness of the model reduction and exactness verification
methods by numerical examples on practical FNNs.  

In this paper, we evaluated the magnitude of signals with respect to the $L_2$ norm.  
However, when dealing with robustness certification problems in deep learning, 
it can be preferable to employ the $L_\infty$-norm.  
We intend to extend the present methods to this setting and compare the resulting performance
with those in \cite{Raghunathan_NIPS2018,Chen_NIPS2020}.  


\end{document}